\documentclass[12pt]{amsart}
\usepackage{amsthm,amsfonts,amssymb,amscd,mathrsfs}

\usepackage[all]{xy}

\newcommand\Aut{\text{Aut}}

\newcommand\Ext{\operatorname{Ext}\nolimits}

\newcommand\Hom{\operatorname{Hom}\nolimits}

\newcommand\im{\operatorname{Im}\nolimits}

\newcommand{\Pic}{\operatorname{Pic}\nolimits}

\newcommand\rk{\text{rk}}
\newcommand\Fl{Fl}
\newcommand{\Sing}{\operatorname{Sing}\nolimits}
\newcommand\Span{\text{Span}}

\newcommand{\PP}{\mathbb{P}}

\newcommand{\GG}{\mathbf{X}}
\newcommand{\EE}{\mathbf{E}}
\newcommand{\LL}{\mathcal{L}}
\newcommand{\OO}{\mathcal{O}}

\renewcommand{\phi}{\varphi}

\newcommand {\ZZ} {\mathbf {Z}}
\newcommand {\CC} {\mathbb {C}}

\newcommand {\NN} {\mathbf {N}}

\newcommand{\refth}[1]{Theorem \ref{#1}}

\newcommand{\refeq}[1]{(\ref{#1})}

\newtheorem{theorem}{Theorem}[section]
\newtheorem{proposition}[theorem]{Proposition}
\newtheorem{prop}[theorem]{Proposition}
\newtheorem{conjecture}[theorem]{Conjecture}
\newtheorem{lemma}[theorem]{Lemma}
\newtheorem{corollary}[theorem]{Corollary}

\theoremstyle{definition}

\newtheorem{proposition-definition}[theorem]{Proposition-Definition}

\begin{document}

\title{Rank 2 vector bundles on ind-grassmannians}

\author[I.Penkov]{\;Ivan~Penkov}

\address{
Jacobs University Bremen\footnote{International University Bremen prior 
to Spring 2007} \\
School of Engineering and Science,
Campus Ring 1,
28759 Bremen, Germany}
\email{i.penkov@iu-bremen.de}

\author[Tikhomirov]{\;Alexander~S.~Tikhomirov}

\address{
Department of Mathematics\\
State Pedagogical University\\
Respublikanskaya Str. 108
\newline 150 000 Yaroslavl, Russia}
\email{astikhomirov@mail.ru}

\begin{flushright}
\begin{tabular}{l}
To Yuri Ivanovich Manin\\
on the occasion of his 70$^{th}$\\
birthday
\end{tabular}

\end{flushright}

\maketitle

\thispagestyle{empty}

\section{Introduction}
\label{Introduction}

The simplest example of an ind-Grassmannian is the infinite projective 
space $\mathbf P^\infty$.
The Barth-Van de Ven-Tyurin (BVT) Theorem, proved more than 30 years 
ago \cite{BV}, \cite{T},
\cite{Sa} (see also a recent proof by A. Coand\u a and G. Trautmann, 
\cite{CT}), claims that any vector bundle of finite rank on $\mathbf 
P^\infty$ is isomorphic to a
direct sum of line bundles. In the last decade natural examples of 
infinite flag varieties (or flag
ind-varieties) have arisen as homogeneous spaces of locally linear 
ind-groups, \cite{DPW},
\cite{DiP}. In the present paper we concentrate our attention to the 
special case of
ind-Grassmannians, i.e. to inductive limits
of Grassmannians of growing dimension. If $V=\displaystyle\bigcup_{n>k} 
V^n$ is a
countable-dimensional vector space, then the ind-variety
$\mathbf G(k;V)=\displaystyle\lim_\to G(k;V^n)$ (or
simply $\mathbf G(k;\infty)$) of $k$-dimensional subspaces of $V$ is of 
course an
ind-Grassmannian: this is the simplest example beyond $\mathbf 
P^\infty=\mathbf G(1;\infty)$.
A significant difference between $\mathbf G(k;V)$ and a general 
ind-Grassmannian
$\mathbf X=\displaystyle\lim_\to G(k_i;V^{n_i})$ defined via a sequence 
of embeddings
\begin{equation}\label{eq1}
G(k_1;V^{n_1})\stackrel{\phi_1}{\longrightarrow}G(k_2;V^{n_2})
\stackrel{\phi_2}{\longrightarrow}\dots\stackrel{\phi_{m-1}}{\longrightarrow}G(k_m;V^{n_m})
\stackrel{\phi_m}{\longrightarrow}\dots,
\end{equation}
is that in general the morphisms $\phi_m$ can have arbitrary
degrees. We say that the ind-Grassmannian $\mathbf X$ is
\emph{twisted} if $\deg\phi_m>1$ for infinitely many $m$, and that
$\mathbf X$ is \emph{linear} if $\deg\phi_m=1$ for almost all $m$.
\begin{conjecture}\label{con1}
Let the ground field be $\CC$, and let $\mathbf E$ be a vector bundle of 
rank $r\in\ZZ_{>0}$ on an
ind-grasmannian $\mathbf X=\displaystyle\lim_\to G(k_m;V^{n_m})$, i.e.
$\mathbf E=\displaystyle\lim_\gets E_m$, where $\{E_m\}$ is an inverse 
system of vector bundles
of (fixed) rank $r$ on $G(k_m;V^{n_m})$. Then
\begin{itemize}
\item[(i)] $\mathbf E$ is semisimple: it is isomorphic to a direct sum 
of simple vector bundles
on $\mathbf X$, i.e. vector bundles on $\mathbf X$ with no non-trivial 
subbundles;
\item[(ii)] for $m\gg0$ the restriction of each simple bundle $\mathbf 
E$ to
$G(k_m,V^{n_m})$ is a homogeneous vector bundle;
\item[(iii)] each simple bundle $\mathbf E'$ has rank 1 unless $\mathbf 
X$ is isomorphic
$\mathbf G(k;\infty)$ for some $k$: in the latter case $\mathbf E'$, 
twisted by a suitable
line bundle,
is isomorphic to a simple subbundle of the tensor algebra 
$T^{\cdot}(\mathbf S)$, $\mathbf S$
being the tautological bundle of rank $k$ on $\mathbf G(k;\infty)$;
\item[(iv)] each simple bundle $\mathbf E$ (and thus each vector bundle 
of finite rank on
$\mathbf X$) is trivial whenever $\mathbf X$ is a twisted 
ind-Grassmannian.
\end{itemize}
\end{conjecture}
The BVT Theorem and Sato's theorem about finite rank bundles on 
$\mathbf G(k;\infty)$,
\cite{Sa}, \cite{Sa2}, as well as the results in \cite{DP}, are 
particular cases of the above conjecture.
The purpose of the present note is to prove Conjecture \ref{con1} for 
vector bundles of rank 2,
and also for vector bundles of arbitrary rank $r$ on linear 
ind-Grassmannians $\mathbf X$.

In the 70's and 80's Yuri Ivanovich Manin taught us mathematics in (and 
beyond) his seminar,
and the theory of vector bundles was a reoccuring topic (among many 
others). In 1980, he asked
one of us (I.P.) to report on A. Tyurin's paper \cite{T}, and most 
importantly to try to
understand this paper. The present note is a very preliminary progress 
report.

\textbf{Acknowledgement. }We acknowledge the support and hospitality of 
the Max Planck Institute for Mathematics in Bonn where the present 
note was conceived. A. S. T. also acknowledges partial support from 
Jacobs University Bremen. Finally, we thank the referee for a number of sharp comments.

\section{Notation and Conventions}
The ground field is $\CC$. Our notation is mostly standard: if $X$ is 
an algebraic variety,
(over $\CC$), $\OO_X$ denotes its structure sheaf, $\Omega^1_X$ 
(respectively $T_X$) denotes
the cotangent (resp. tangent) sheaf on X under the assumption that $X$ 
is smooth etc.
If $F$ is a sheaf on $X$, its cohomologies are denoted by $H^i( F)$, 
$h^i(F):=\dim H^i(F)$, and
$\chi(F)$ stands for the Euler characteristic of $F$.
The Chern classes of $F$ are denoted by $c_i(F)$. If $f:X\to Y$ is a 
morphism, $f^*$ and $f_*$
denote respectively the inverse and direct image functors of 
$\OO$-modules. All vector bundles
are assumed to have finite rank. We denote the dual of a sheaf of 
$\mathcal O_X$-modules $F$
(or that of a vector space) by the superscript $^\vee$.
Furthermore, in what follows for any ind-Grassmannian $\mathbf X$ 
defined by \refeq{eq1}, no
embedding $\phi_i$ is an isomorphism.

We fix a finite dimensional space $V$ and denote by $X$ the 
Grassmannian $G(k;V)$ for $k<\dim V$. In the sequel we write sometimes $G(k;n)$ 
indicating simply the dimension of $V$. Below we will often consider 
(parts of) the following diagram of flag varieties:
\begin{equation}\label{eqDiag}
\xymatrix{
&&Z:=\Fl(k-1,k,k+1;V) \ar[ld]_{\pi_1} \ar[dr]^{\pi_2}  & \\
&Y:=\Fl(k-1,k+1;V)\ar[ld]_{p_1}\ar[rd]^{p_2}&&X:=G(k;V), \\
Y^1:=G(k-1;V)&&Y^2:=G(k+1;V)&\\
}
\end{equation}
under the assumption that $k+1<\dim V$. Moreover we reserve the letters 
$X,Y,Z$ for the
varieties in the above diagram. By $S_k$, $S_{k-1}$, $S_{k+1}$ we 
denote the tautological
bundles on $X$,$Y$ and $Z$, whenever they are defined ($S_k$ is defined 
on
$X$ and $Z$, $S_{k-1}$ is defined on $Y^1$, $Y$ and $Z$, etc.). By
$\mathcal O_X(i)$, $i\in \ZZ$, we denote the isomorphism class (in the 
Picard group $\Pic X$)
of the line bundle $(\Lambda^k(S_k^\vee))^{\otimes i}$, where 
$\Lambda^k$ stands for the
$k^{th}$ exterior power (in this case maximal exterior power as $\rk 
S_k^\vee=k$). The Picard
group of $Y$ is isomorphic to the direct product of the Picard groups 
of $Y^1$ and $Y^2$,
and by $\OO_Y(i,j)$ we denote the isomorphism class of the line bundle
$p_1^*(\Lambda^{k-1}(S_{k-1}^\vee))^{\otimes i}
\otimes_{\OO_Y}p_2^*(\Lambda^{k+1}(S_{k+1}^\vee))^{\otimes j}$.

If $\phi:X=G(k;V)\to X':=G(k;V')$ is an embedding, then 
$\phi^*\OO_{X'}(1)\simeq \OO_X(d)$
for some $d\in\ZZ_{\geq 0}$: by definition $d$ is the \emph{degree} 
$\deg\phi$ of $\phi$.
We say that $\phi$ is linear if $\deg\phi=1$. By a \textit{projective 
subspace}
(in particular a \emph{line}, i.e. a 1-dimensional projective subspace) 
of $X$ we mean a
linearly embedded projective space into $X$. It is well known that all 
such are Schubert
varieties of the form
$\{V^k\in X| V^{k-1}\subset V^k\subset V^t\}$ or $\{V^k\in X| 
V^i\subset V^k\subset V^{k+1}\}$,
where $V^k$ is a variable $k$-dimensional subspace of $V$, and
$V^{k-1}$, $V^{k+1}$, $V^t$, $V^i$ are fixed subspaces of $V$ of 
respective dimensions
$k-1$, $k+1$, $t$, $i$. (Here and in what follows $V^t$ always denotes 
a vector space of
dimension $t$). In other words, all projective subspaces of $X$ are of 
the form
$G(1;V^t/V^{k-1})$ or $G(k-i, V^{k+1}/V^i)$.
Note also that $Y=\Fl(k-1,k+1;V)$ is the variety of lines in 
$X=G(k;V)$.

\section{The linear case}
We consider the cases of linear and twisted ind-Grassmannians 
separately. In the case of a
linear ind-Grassmannian, we show that Conjecture \ref{con1} is a 
straightforward corollary of
existing results combined with the following proposition. We recall, 
\cite{DP}, that a
\textit{standard extension} of Grassmannians is an embedding of the 
form
\begin{equation}\label{eq31}
G(k;V)\to G(k+a;V\oplus \hat W), \quad \{ V^k\subset 
\CC^n\}\mapsto\{V^k\oplus W\subset V\oplus\hat W\},
\end{equation}
where $W$ is a fixed $a$-dimensional subspace of a finite dimensional 
vector space $\hat W$.

\begin{proposition}\label{linear embed}
Let $\phi:X=G(k;V)\to X':=G(k';V')$ be an embedding of degree 1. Then 
$\phi$ is a standard extension, or $\phi$ factors through a standard 
extension $\PP^r\to G(k';V')$ for some $r$.
\end{proposition}
\begin{proof}
We assume that $k\leq n-k$, $k\leq n'-k'$, where $n=\dim V$ and 
$n'=\dim V'$, and use induction on $k$. For $k=1$ the statement is obvious as 
the image of $\phi$ is a projective subspace of $G(k';V')$ and hence 
$\phi$ is a standard extension. Assume that the statement is true for 
$k-1$. Since $\deg \phi=1$, $\phi$ induces an embedding $\phi_Y:Y\to Y'$, 
where $Y=\Fl(k-1,k+1;V)$ is the variety of lines in $X$ and 
$Y\:=\Fl(k'-1,k'+1;V')$ is the variety of lines in $X'$. Moreover, clearly we have 
a commutative diagram of natural projections and embeddings
\[
  \xymatrix{
&Z\ar[rrr]^{\phi_Z}\ar[dl]_{\pi_1}\ar[dr]^{\pi_2}&&&Z'\ar[dl]_{\pi_1'}\ar[dr]^{\pi_2'}& 
\\
Y\ar[dr]&&X\ar[dr]&Y'&&X',\\
&\ar[r]_{\phi_Y}&\ar[ur]&\ar[r]_{\phi}&\ar[ur]&
}
\]
where $Z:=\Fl(k-1,k,k+1;V)$ and $Z':=\Fl(k'-1,k',k'+1;V')$.

We claim that there is an isomorphism
\begin{equation}\label{eqLE1}
\phi^*_Y\OO_{Y'}(1,1)\simeq\OO_Y(1,1).
\end{equation}
Indeed, $\phi^*_Y\OO_{Y'}(1,1)$ is determined up to isomorphism by its 
restriction to the fibers of $p_1$ and $p_2$ (see diagram 
\refeq{eqDiag}), and therefore it is enough to check that
\begin{equation}\label{eqLE2}
\phi^*_Y\OO_{Y'}(1,1)_{|p_1^{-1}(V^{k-1})}\simeq\OO_{p_1^{-1}(V^{k-1})}(1),
\end{equation}
\begin{equation}\label{eqLE21}
\phi^*_Y\OO_{Y'}(1,1)_{|p_2^{-1}(V^{k+1})}\simeq 
\OO_{p_2^{-1}(V^{k+1})}(1)
\end{equation}
for some fixed subspaces $V^{k-1}\subset V$, $V^{k+1}\subset V$. Note 
that
the restriction of $\phi$ to the projective subspace 
$G(1;V/V^{k-1})\subset X$
is simply an isomorphism of $G(1;V/V^{k-1})$ with a projective subspace 
of $X'$,
hence the map induced by $\phi$ on the variety $G(2;V/V^{k-1})$ of 
projective
lines in $G(1;V/V^{k-1})$ is an isomorphism with the Grassmannian of 
2-dimensional
subspaces of an appropriate subquotient of $V'$. Note furthermore 
that $p_1^{-1}(V^{k-1})$ is nothing but the variety of lines 
$G(2;V/V^{k-1})$ in $G(1;V/V^{k-1})$, and that the image of $G(2;V/V^{k-1})$ under 
$\phi$ is nothing but $\phi_Y(p_1^{-1}(V^{k-1}))$. This shows that  the 
restriction of $\phi^*_Y\OO_{Y'}(1,1)$ to $G(2;V/V^{k-1})$ is 
isomorphic to the restriction of $\OO_Y(1,1)$ to $G(2;V/V^{k-1})$, and we obtain 
\refeq{eqLE2}. The isomorphism \refeq{eqLE21} follows from a very 
similar argument.

The isomorphism \refeq{eqLE1} leaves us with two alternatives:
\begin{equation}\label{eqLE3}
\phi^*_{Y}\OO_{Y'}(1,0)\simeq\OO_Y \mathrm{~or~} 
\phi_Y^*\OO_{Y'}(0,1)\simeq \OO_Y,
\end{equation}
or
\begin{equation}\label{eqLE4}
\phi^*_{Y}\OO_{Y'}(1,0)\simeq\OO_Y(1,0) \mathrm{~or~} 
\phi_Y^*\OO_{Y'}(1,0)\simeq \OO_Y(0,1).
\end{equation}
Let \refeq{eqLE3} hold, more precisely let 
$\phi_Y^*\OO_{Y'}(1,0)\simeq\OO_Y$.
Then $\phi_Y$ maps each fiber of $p_2$ into a single point in $Y'$ 
(depending on the image in
$Y^2$ of this fiber), say $({(V')}^{k'-1}\subset {(V')}^{k'+1})$, and 
moreover the space
${(V')}^{k'-1}$ is constant. Thus $\phi$ maps $X$ into the projective 
subspace $G(1;V'/{(V')}^{k'-1})$ of $X'$. If 
$\phi_Y^*\OO_{Y'}(0,1)\simeq\OO_Y$, then $\phi$ maps $X$ into the projective subspace 
$G(1;{(V')}^{k'+1})$ of $X'$. Therefore, the Proposition is proved in the case 
\refeq{eqLE3} holds.

We assume now that \refeq{eqLE4} holds. It is easy to see that 
\refeq{eqLE4} implies that
$\phi$ induces a linear embedding $\phi_{Y^1}$ of $Y^1:=G(k-1;V)$ into 
$G(k'-1;V')$ or
$G(k'+1;V')$. Assume that $\phi_{Y^1}:Y^1\to {(Y')}^1:=G(k'-1;V')$ (the 
other case is
completely similar). Then, by the induction assumption, $\phi_{Y^1}$ is 
a standard extension
or factors through a standard extension $\PP^r\to {(Y')}^1$. If 
$\phi_{Y^1}$ is a standard
extension corresponding to a fixed subspace $W\subset \hat W$, then 
$\phi_{Y^1}^* S_{k'-1}\simeq S_{k-1}\oplus 
\left(W\otimes_\CC\OO_{Y^1}\right)$ and we have a vector bundle monomorphism
\begin{equation}\label{eqLE5}
0\to\pi_1^*p_1^*\phi_{Y^1}^*S_{k'-1}\to \pi_2^*\phi^*S_{k'}.
\end{equation}
By restricting \refeq{eqLE5} to the fibers of $\pi_1$ we see that the 
quotient line bundle 
$\pi_2^*\phi^*S_{k'}/\pi_1^*p_1^*\phi_{Y^1}^*S_{k'-1}$ is isomorphic to $S_k/S_{k-1}\otimes \pi_1^*p_1^*\LL$, where $\LL$ 
is a line bundle on $Y^1$. Applying $\pi_{2*}$ we obtain
\begin{equation}\label{eqLE6}
0\to W\otimes_\CC \OO_X\to\pi_{2*}(\pi_2^*\phi^*S_{k'})=\phi^*S_{k'}\to 
\pi_{2*}((S_k/S_{k-1})\otimes\pi_1^*p_1^*\LL) \to 0.
\end{equation}
Since $\rk\phi^*S_{k'}=k'$ and $\dim W=k'-k$, 
$\rk\pi_{2*}((S_k/S_{k-1})\otimes\pi_1^*p_1^*\LL)=k$, which implies immediately that $\LL$ is 
trivial. Hence \refeq{eqLE6} reduces to  $0\to 
W\otimes_{\CC}\OO_X\to\phi^*S_{k'}\to S_k\to 0$, and thus
\begin{equation}\label{eqLE7}
\phi^*S_{k'}\simeq S_k\oplus \left(W\otimes_\CC\OO_X\right)
\end{equation}
as there are no non-trivial extensions of $S_k$ by a trivial bundle. 
Now \refeq{eqLE7} implies
that $\phi$ is a standard extension.

It remains to consider the case when $\phi_{Y^1}$ maps $Y^1$ into a 
projective subspace
$\PP^s$ of ${(Y')}^1$. Then $\PP^s$ is of the form 
$G(1;V'/{(V')}^{k'-2})$ for some
${(V')}^{k'-2}\subset V'$, or of the form $G(k'-1;{(V')}^{k'})$ for 
some
${(V')}^{k'}\subset V'$. The second case is clearly impossible because 
it would imply that
$\phi$ maps $X$ into the single point ${(V')}^{k'}$. Hence 
$\PP^s=G(1;V'/{(V')}^{k'-2})$ and
$\phi$ maps $X$ into the Grassmannian $G(2;V'/{(V')}^{k'-2})$ in 
$G(k';V')$. Let $S_2'$ be the
rank 2 tautological bundle on $G(2;V'/{(V')}^{k'-2})$. Then its 
restriction
$S'':=\phi^*S_2'$ to any line $l$ in $X$ is isomorphic to 
$\OO_{l}\oplus\OO_{l}(-1)$, and we
claim that this implies one of the two alternatives:
\begin{equation}\label{eqLE8}
S''\simeq\OO_X\oplus\OO_X(-1)
\end{equation}
or
\begin{equation}\label{eqLE9}
S''\simeq S_2 \text{~and~} k=2,\text{~or~} S''\simeq(V\otimes_\CC 
\OO_X)/S_2\text{~and~}k=n-k=2.
\end{equation}
Let $k\geq 2$. The evaluation map $\pi_1^*\pi_{1*}\pi_2^*S''\to 
\pi_2^*S''$ is a monomorphism
of the line bundle $ \pi_1^*\LL:=\pi_1^*\pi_{1*}\pi_2^*S''$ into 
$\pi_2^*S''$
(here $\LL:=\pi_{1*}\pi_2^*S''$). Restricting this monomorphism to the 
fibers of $\pi_2$ we see
immediately that $\pi_1^*\LL$ is trivial when restricted to those 
fibers and is hence trivial.
Therefore $\LL$ is trivial, i.e. $\pi_1^*\LL=\OO_Z$. Push-down to $X$ 
yields
\begin{equation}\label{eqLE10}
0\to\OO_X\to S''\to\OO_X(-1)\to 0,
\end{equation}
and hence \refeq{eqLE10} splits as $\Ext^1(\OO_X(-1),\OO_X)=0$. 
Therefore \refeq{eqLE8} holds.
For $k=2$, there is an additional possibility for the above 
monomorphisms to be of the form
$\pi_1^*\OO_Y(-1,0)\to\pi_2^*S$ (or of the form 
$\pi_1^*\OO_Y(0,-1)\to\pi_2^*S$ if $n-k=2$)
which yields the option \refeq{eqLE9}.

If \refeq{eqLE8} holds, $\phi$ maps $X$ into an appropriate projective 
subspace of $G(2;V'/{(V')}^{k'-2})$ which is then a projective subspace 
of $X'$, and if \refeq{eqLE9} holds, $\phi$ is a standard extension 
corresponding to a zero dimensional space $W$. The proof is now complete.
\end{proof}

We are ready now to prove the following theorem.
\begin{theorem} Conjecture \ref{con1} holds for any linear 
ind-Grassmannian $\mathbf X$.
\end{theorem}
\begin{proof}
Assume that $\deg \phi_m=1$ for all $m$, and apply Proposition 
\ref{linear embed}. If
infinitely many $\phi_m$'s factor through respective projective 
subspaces, then $\mathbf X$ is
isomorphic to $\mathbf P^\infty$ and the BVT Theorem implies Conjecture 
\ref{con1}. Otherwise,
all $\phi_m$'s are standard extensions of the form \refeq{eq31}. There 
are two alternatives:
$\displaystyle\lim_{m\to\infty} 
k_{m}=\lim_{m\to\infty}(n_{m}-k_{m})=\infty$, or one of the
limits $\displaystyle\lim_{m\to \infty}k_{m}$ or 
$\displaystyle\lim_{m\to \infty}(n_{m}-k_{m})$
equals $l$ for some $l\in \NN$. In the first case the claim of 
Conjecture \ref{con1} is proved
in \cite{DP}: Theorem 4.2. In the second case $\mathbf X$ is isomorphic 
to
$\mathbf G(l;\infty)$, and therefore Conjecture \ref{con1} is proved in 
this case by E. Sato
in \cite{Sa2}.
\end{proof}

\section{Auxiliary results}

In order to prove Conjecture \ref{con1} for rank 2 bundles $\mathbf E$ 
on a twisted
ind-Grassmannian $\mathbf X=\displaystyle \lim_\to G(k_m;V^{n_m})$, we 
need to prove that the
vector bundle $\mathbf E=\displaystyle\lim_{\gets}E_m$ of rank 2 on 
$\mathbf X$ is trivial,
i.e. that $E_m$ is a trivial bundle on $G(k_m;V^{n_m})$ for each $m$.
From this point on we assume that none of the Grassmannians 
$G(k_m;V^{n_m})$ is a projective
space, as for a twisted projective ind-space Conjecture 1.1 is proved 
in \cite{DP} for bundles
of arbitrary rank $r$.

The following known proposition gives a useful triviality criterion for 
vector bundles of
arbitrary rank on Grassmannians.

\begin{prop}\label{prop31}
A vector bundle $E$ on $X=G(k;n)$ is trivial iff its restriction 
$E_{|l}$ is trivial for every
line $l$ in $G(k;n)$, $l\in Y=\Fl(k-1,k+1;n)$.
\end{prop}
\begin{proof}
We recall the proof given in \cite{P}. It uses the well known fact that 
the Proposition holds
for any projective space, [OSS, Theorem 3.2.1]. Let first $k=2$, $n=4$, 
i.e. $X=G(2;4)$. Since
$E$ is linearly trivial, $\pi_2^*E$ is trivial along the fibers of 
$\pi_1$ (we refer here to
diagram \refeq{eqDiag}). Moreover, $\pi_{1*}\pi_2^*E$ is trivial along 
the images of the fibers
of $\pi_2$ in $Y$. These images are of the form $\PP_1^1\times\PP_2^1$, 
where $\PP_1^1$
(respectively $\PP_2^1$) are lines in $Y^1:=G(1;4)$ and $Y^2:=G(3;4)$. 
The fiber of $p_1$ is
filled by lines of the form $\PP^1_2$, and thus $\pi_{1*}\pi_2^*E$ is 
linearly trivial, and
hence trivial along the fibers of $p_1$. Finally the lines of the form 
$\PP_1^1$ fill $Y^1$,
hence ${p_1}_*\pi_{1*}\pi_2^*E$ is also a trivial bundle. This implies 
that
$E=\pi_{2*}\pi_1^*p_1^*(p_{1*}\pi_{1*}\pi_2^*E)$ is also trivial.

The next case is the case when $k=2$ and $n$ is arbitrary, $n\geq 5$. 
Then the above argument
goes through by induction on $n$ since the fiber of $p_1$ is isomorphic 
to $G(2;n-1)$. The
proof is completed by induction on $k$ for $k\geq 3$: the base of $p_1$ 
is $G(k-1;n)$ and the
fiber of $p_1$ is $G(2;n-1)$.
\end{proof}

If $C\subset N$ is a smooth rational curve in an algebraic variety $N$ 
and $E$ is a vector
bundle on $N$, then by a classical theorem of Grothendieck, 
$\displaystyle E_{|C}$ is isomorphic
to $\bigoplus_i\OO_C(d_i)$ for some $d_1\geq d_2\geq\dots\geq d_{\rk 
E}$. We call the ordered
$\rk E$-tuple $(d_1,\dots,d_{\rk E})$ \emph{the splitting type} of 
$E_{|C}$ and denote it by
$\mathbf{d}_E(C)$. If $N=X=G(k;n)$, then the lines on $N$ are 
parametrized by points $l\in Y$,
and we obtain a map
\[
Y\to \ZZ^{\rk E}\ :\ l\mapsto \mathbf{d}_E(l).
\]
By semicontinuity (cf. \cite[Ch.I, Lemma 3.2.2]{OSS}), there is a dense 
open set $U_E\subset Y$ of lines with minimal splitting type with 
respect to the lexicographical ordering on $\ZZ^{\rk E}$. Denote this 
minimal splitting type by $\mathbf{d}_E$. By definition, $U_E=\{l\in Y|~ 
\mathbf{d}_E(l)=\mathbf{d}_E\}$ is the set of \emph{non-jumping} lines of 
$E$, and its complement $Y\setminus U_E$ is the proper closed set of 
\emph{jumping} lines.

A coherent sheaf $F$ over a smooth irreducible variety $N$ is called 
$normal$ if for every open
set $U\subset N$ and every closed algebraic subset $A\subset U$ of 
codimension at least 2 the
restriction map ${F}(U)\to {F}(U\smallsetminus A)$ is surjective. It is 
well known that, since
$N$ is smooth, hence normal, a normal torsion-free sheaf $F$ on $N$ is 
reflexive, i.e.
$F^{\lor\lor}=F$. Therefore, by \cite[Ch.II, Theorem 2.1.4]{OSS} $F$ is 
necessarily a line
bundle (see \cite[Ch.II, 1.1.12 and 1.1.15]{OSS}).

\begin{theorem}\label{thSubbdl}
Let $E$ be a rank $r$ vector bundle of splitting type
$\mathbf{d}_E=(d_1,...,d_r),\ d_1\ge...\ge d_r,$ on $X=G(k;n)$.
If $d_s-d_{s+1}\ge2$ for some $s<r$, then there is a normal subsheaf 
$F\subset E$ of rank $s$
with the following properties: over the open set 
$\pi_2(\pi_1^{-1}(U_E))\subset X$ the
sheaf $F$ is a subbundle of $E$, and for any $l\in U_E$
$$
F_{|l}\simeq\overset{s}{\underset{i=1}\bigoplus}\mathcal{O}_{l}(d_i).
$$
\end{theorem}
\begin{proof}
It is similar to the proof of Theorem 2.1.4 of \cite[Ch.II]{OSS}. 
Consider the vector bundle
$E'=E\bigotimes\mathcal{O}_X(-d_s)$ and the evaluation map
$\Phi:\pi_1^*\pi_{1*}\pi_2^*E'\to \pi_2^*E'$. The definition of $U_E$ 
implies that
$\Phi_{|\pi_1^{-1}(U_E)}$ is a morphism of constant rank $s$ and that 
its image
${\rm \im}\Phi\subset \pi_2^*E'$ is a subbundle of rank $s$ over 
$\pi_1^{-1}(U_E)$.
Let $M:=\pi_2^*E'/{\rm im}\Phi$, let $T(M)$ be the torsion subsheaf of 
$M$, and
$F':=\ker(\pi_2^*E'\to M':=M/T(M))$. Consider the singular set $\Sing 
F'$ of the
sheaf $F'$ and set $A:=Z\smallsetminus\Sing F'$. By the above, $A$ is 
an open subset
of $Z$ containing $\pi_1^{-1}(U_E)$ and $f={\pi_2}_{|A}:A\to 
B:=\pi_2(A)$ is a submersion
with connected fibers.

Next, take any point $l\in Y$ and put $ L:=\pi_1^{-1}(l)$. By 
definition,
$L\simeq\mathbb{P}^1$, and we have
\begin{equation}\label{tangent}
{T_{Z/X}}_{|L}\simeq\mathcal{O}_{L}(-1)^{\oplus(n-2)},
\end{equation}
where $T_{Z/X}$ is the relative tangent bundle of Z over X. The 
construction of the sheaves
$F'$ and $M$ implies that for any
$l\in U_E$: 
${F'}^{\vee}_{|{L}}=\oplus_{i=1}^s\mathcal{O}_{L}(-d_i+d_s),\ \
{M'}_{|{L}} =\oplus_{i=s+1}^r\mathcal{O}_{L}(d_i-d_s)$.
This, together with (\ref{tangent}) and the condition 
$d_s-d_{s+1}\ge2,$ immediately implies
that $H^0(\Omega^1_{A/B}\otimes{F'}^{\vee}\otimes M'_{|{L}})=0$. Hence
$H^0(\Omega^1_{A/B}\otimes{F'}^{\vee}\otimes M'_{|\pi_1^{-1}(U_E)})=0$, 
and thus, since
$\pi_1^{-1}(U_E)$ is dense open in $Z$,
$\Hom(T_{A/B},\mathcal H om(F',M'_{|A}))=
H^0(\Omega^1_{A/B}\otimes{F'}^{\vee}\otimes M'_{|A})=0.$
Now we apply the Descent Lemma (see \cite[Ch.II, Lemma 2.1.3]{OSS}) to 
the data
$(f_{|\pi_1^{-1}(U_E)}:\pi_1^{-1}(U_E)\to V_E,\ F'_{|\pi_1^{-1}(U_E)}
\subset E'_{|\pi_1^{-1}(U_E)})$. Then
$F:=(\pi_{2*}F')\otimes\mathcal{O}_X(-d_s)$
is the desired sheaf.
\end{proof}

\section{The case $\rk\EE=2$}
In what follows, when considering a twisted ind-Grassmannian $\mathbf 
X=\displaystyle\lim_\to G(k_m;V^{n_m})$ we set $G(k_m;V^{n_m})=X_m$.  
\refth{thSubbdl} yields now the following corollary.
\begin{corollary}\label{d=(0,0)}
Let $\displaystyle\mathbf{E}=\lim_{\gets}E_m$ be a rank 2 vector bundle 
on a twisted
ind-Grassmannian $\displaystyle\mathbf{X}=\lim_{\to}X_m$. Then there 
exists $m_0\ge1$ such that $\mathbf{d}_{E_m}=(0,0)$ for any $m\ge m_0.$
\end{corollary}
\begin{proof}
Note first that the fact that $\mathbf X$ is twisted implies
\begin{equation}\label{c_1=0}
c_1(E_m)=0,\ m\ge1.
\end{equation}
Indeed, $c_1(E_m)$ is nothing but the integer corresponding to the line 
bundle $\Lambda^2(E_m)$ in the identification of $\Pic X_m$ with $\ZZ$. 
As $\mathbf X$ is twisted, 
$c_1(E_m)=\deg\phi_m\deg\phi_{m+1}\dots\deg\phi_{m+k}c_1(E_{m+k+1})$ for any $k\geq 1$, in other words $c_1(E_m)$ 
is divisible by larger and larger integers and hence $c_1(E_m)=0$ (cf. 
\cite[Lemma 3.2]{DP}). Suppose that for any $m_0\ge1$ there exists $m\ge 
m_0$ such that $\mathbf{d}_{E_m}=(a_m,-a_m)$ with $a_m>0$. Then Theorem 
\ref{thSubbdl} applies to $E_m$ with $s=1$, and hence $E_m$ has a 
normal rank-1 subsheaf $F_m$ such that
\begin{equation}\label{F|l}
F_{m|l}\simeq\mathcal{O}_{l}(a_m)
\end{equation}
for a certain line $l$ in $X_m$. Since $F_m$ is a torsion-free normal 
subsheaf of the vector bundle $E$, the sheaf $F_m$ is a line bundle, 
i.e. $F_m\simeq\OO_{X_m}(a_m)$. Therefore we have a monomorphism:
\begin{equation}\label{injectn}
0\to\mathcal{O}_{X_m}(a_m)\to E_m,\ \ \  a_m\ge1.
\end{equation}
This is clearly impossible. In fact, this monomorphism implies in view 
of (\ref{c_1=0}) that any rational curve $C\subset X_m$ of degree 
$\delta_m:=\deg\phi_1\cdot...\cdot\deg\phi_{m-1}$ has splitting type 
$\mathbf{d}_{E_m}(C)=(a'_m,-a'_m)$, where $a'_m\ge a_m\delta_m\ge\delta_m$. 
Hence, by semiconinuity, any line $l\in X_1$ has splitting type 
$\mathbf{d}_{E_1}(l)=(b,-b),\ \ b\ge\delta_m$. Since $\delta_m\to\infty$ as 
$m_0\to\infty,$ this is a contradiction.
\end{proof}

We now recall some standard facts about the Chow rings of 
$X_m=G(k_m;V^{n_m}),$
(see, e.g., \cite[14.7]{F}):
\begin{itemize}
\item[(i)] $A^1(X_m)=\Pic(X_m)=\mathbb{Z}[\mathbb{V}_m]$,
$A^2(X_m)=\mathbb{Z}[\mathbb{W}_{1,m}]\oplus\mathbb{Z}[\mathbb{W}_{2,m}]$, 
where
$\mathbb{\mathbb{V}}_m,\mathbb{W}_{1,m},\mathbb{W}_{2,m}$ are the 
following
Schubert varieties:
$\mathbb{V}_m:=\{V^{k_m}\in X_m|\ \dim(V^{k_m}\cap V_0^{n_m-k_m})\ge1$
for a fixed subspace
$V_0^{n_m-k_m-1}$ of $V^{n_m}\}$,
$\mathbb{W}_{1,m}:=\{V^{k_m}\in X_m| $ $\dim (V^{k_m}\cap 
V_0^{n_m-k_m-1})\ge1$
for a fixed subspace $V_0^{n_m-k_m-1}$ in $V^{n_m}\}$,
$\mathbb{W}_{2,m}:=\{{V}^{k_m}\in X_m|\ \dim({V}^{k_m}\cap 
V_0^{n_m-k_m+1})\ge2$
for a fixed subspace $V_0^{n_m-k_m+1}$ of $V^{n_m}\}$;
\item[(ii)] $[\mathbb{V}_m]^2=[\mathbb{W}_{1,m}]+[\mathbb{W}_{2,m}]$ in 
$A^2(X_m)$;
\item[(iii)] 
$A_2(X_m)=\mathbb{Z}[\mathbb{P}^2_{1,m}]\oplus\mathbb{Z}[\mathbb{P}^2_{2,m}]$,
where the projective planes
$\mathbb{P}^2_{1,m}$ (called \emph{$\alpha$-planes}) and 
$\mathbb{P}^2_{2,m}$
(called \emph{$\beta$-planes}) are respectively the Schubert varieties
$\mathbb{P}^2_{1,m}:=\{V^{k_m}\in X_m|\ V_0^{k_m-1}\subset 
{V}^{k_m}\subset V_0^{k_m+2}$
for a fixed flag $V_0^{k_m-1}\subset V_0^{k_m+2}$ in $V^{n_m}\}$, 
$\mathbb{P}^2_{2,m}:=\{V^{k_m}\in X_m|\ V_0^{k_m-2}\subset {V}^{k_m}\subset 
V_0^{k_m+1}$ for a fixed flag $V_0^{k_m-2}\subset V_0^{k_m+1}$ in 
$V^{n_m}\};$
\item[(iv)] the bases $[\mathbb{W}_{i,m}]$ and $[\mathbb{P}^2_{j,m}]$
are dual in the standard sense that
$[\mathbb{W}_{i,m}]\cdot[\mathbb{P}^2_{j,m}]=\delta_{i,j}.$
\end{itemize}
\begin{lemma}\label{c_2(E_m)=0}
There exists $m_1\in\ZZ_{>0}$ such that for any $m\ge m_1$ one of the 
following holds:
\begin{itemize}
\item[(1)] $c_2({E_m}_{|\mathbb{P}^2_{1,m}})>0,$ 
$c_2({E_m}_{|\mathbb{P}^2_{2,m}})\le0$,
\item[(2)] $c_2({E_m}_{|\mathbb{P}^2_{2,m}})>0,$ 
$c_2({E_m}_{|\mathbb{P}^2_{1,m}})\le0$,
\item[(3)] $c_2({E_m}_{|\mathbb{P}^2_{1,m}})=0$, 
$c_2({E_m}_{|\mathbb{P}^2_{2,m}})=0$.
\end{itemize}
\end{lemma}

\begin{proof}
According to (i), for any $m\ge1$ there exist 
$\lambda_{1m},\lambda_{2m}\in\ZZ$ such that
\begin{equation}\label{c_2(E_m)}
c_2(E_m)=\lambda_{1m}[\mathbb{W}_{1,m}]+\lambda_{2m}[\mathbb{W}_{2,m}].
\end{equation}
Moreover, (iv) implies
\begin{equation}\label{lambda_jm}
\lambda_{jm}=c_2({E_m}_{|\mathbb{P}^2_{j,m}}),\ \ j=1,2.
\end{equation}
Next, (i) yields:
\begin{equation}\label{abcd}
\phi_m^*[\mathbb{W}_{1,m+1}]=a_{11}(m)[\mathbb{W}_{1,m}]+a_{21}(m)[\mathbb{W}_{2,m}],\ 
\
\phi_m^*[\mathbb{W}_{2,m+1}]=a_{12}(m)[\mathbb{W}_{1,m}]+a_{22}(m)[\mathbb{W}_{2,m}],
\end{equation}
where $a_{ij}(m)\in\mathbb{Z}$. Consider the $2\times2$-matrix 
$A(m)=(a_{ij}(m))$
and the column vector $\Lambda_m=(\lambda_{1m},\lambda_{2m})^t.$ Then, 
in view of (iv),
the relation (\ref{abcd}) gives: $\Lambda_m=A(m)\Lambda_{m+1}$.
Iterating this equation and denoting by $A(m,i)$ the $2\times2$-matrix
$A(m)\cdot A(m+1)\cdot...\cdot A(m+i),\ i\ge1,$ we obtain
\begin{equation}\label{Lambda_m}
\Lambda_m=A(m,i)\Lambda_{m+i+1}.
\end{equation}
The twisting condition 
$\phi_m^*[\mathbb{V}_{m+1}]=\deg\phi_m[\mathbb{V}_{m}]$ together with (ii) implies:
$\phi_m^*([\mathbb{W}_{1,m+1}]+[\mathbb{W}_{2,m+1}])=(\deg\phi_m)^2([\mathbb{W}_{1,m}]+[\mathbb{W}_{2,m}])$.
Substituting (\ref{abcd}) into the last equality, we have:
$a_{11}(m)+a_{12}(m)=a_{21}(m)+a_{22}(m)=(\deg\phi_m)^2,\ \ \ m\ge1.$
This means that the column vector ${v}=(1,1)^t$ is an eigenvector of 
$A(m)$ with eigenvalue $(\deg\phi_m)^2$. Hence, it is an 
eigenvector of $A(m,i)$ with the
eigenvalue
$d_{m,i}=(\deg\phi_m)^2(\deg\phi_{m+1})^2...(\deg\phi_{m+i})^2:$
\begin{equation}\label{eigen}
A(m,i){v}=d_{m,i}{v}.
\end{equation}
Notice that the entries of $A(m),\ m\ge1,$ are nonnegative integers
(in fact, from the definition of the Schubert varieties 
$\mathbb{W}_{j,m+1}$ it immediately follows
that $\phi_m^*[\mathbb{W}_{j,m+1}]$ is an effective cycle on $X_m$, so 
that (\ref{abcd}) and
(iv) give 
$0\le\phi_m^*[\mathbb{W}_{i,m+1}]\cdot[\mathbb{P}^2_{j,m}]=a_{ij}(m)$);
hence also the entries of $A(m,i),\ m,i\ge1,$ are nonnegative 
integers).
Besides, clearly $d_{m,i}\to\infty$ as $i\to\infty$ for any $m\ge1$.
This, together with (\ref{Lambda_m}) and (\ref{eigen}), implies that,
for $m\gg1$, $\lambda_{1m}$ and $\lambda_{2m}$ cannot both be nonzero 
and have the same sign.
This together with (\ref{lambda_jm}) is equivalent to the statement of 
the Lemma.
\end{proof}

In what follows we denote the $\alpha$-planes and the $\beta$-planes on 
$X=G(2;4)$ respectively by $\PP_\alpha^2$ and $\PP_\beta^2$.
\begin{proposition}\label{not exist}
There exists no rank 2 vector bundle $E$ on the Grassmannian $X=G(2;4)$ 
such that:
\begin{itemize}
\item[(a)] $c_2(E)=a[\mathbb{P}^2_{\alpha}],\ \ a>0,$
\item[(b)] $E_{|\mathbb{P}^2_{\beta}}$ is trivial for a generic 
$\beta$-plane
$\mathbb{P}^2_{\beta}$ on $X$.
\end{itemize}
\end{proposition}

\begin{proof} Now assume that there exists a vector bundle $E$ on $X$ 
satisfying the conditions (a) and (b)
of the Proposition. Fix a $\beta$-plane $P\subset X$ such that
\begin{equation}\label{E|Y}
E_{|P}\simeq\mathcal{O}_{P}^{\oplus2}.
\end{equation}
As $X$ is the Grassmannian of lines in $\mathbb{P}^3$, the plane $P$ is 
the dual plane of a
certain plane $\tilde P$ in $\mathbb{P}^3$. Next, fix a point
$x_0\in\mathbb{P}^3\smallsetminus\tilde P$ and denote by $S$ the 
variety of lines in $\PP^3$  which contain $x_0$. Consider the variety
$Q=\{(x,l)\in\mathbb{P}^3\times X\ |\ x\in l\cap\tilde P\}$ with 
natural projections
$p:Q\to S:(x,l)\mapsto\Span(x,x_0)$ and $\sigma:Q\to X:(x,l)\mapsto l$. 
Clearly, $\sigma$ is
the blowing up of $X$ at the plane $P$, and the exceptional divisor 
$D_P=\sigma^{-1}(P)$ is
isomorphic to the incidence subvariety of $P\times\tilde{P}$. Moreover, 
one easily checks that
$Q\simeq\mathbb{P}(\mathcal{O}_{S}(1)\oplus T_{S}(-1))$, so that the 
projection $p:Q\to S$
coincides with the structure morphism 
$\mathbb{P}(\mathcal{O}_{S}(1)\oplus T_{S}(-1))\to S$.
Let $\mathcal{O}_Q(1)$ be the Grothendieck line bundle on $Q$ such
that $p_*\mathcal{O}_Q(1)=\mathcal{O}_{S}(1)\oplus T_{S}(-1)$.
Using the Euler exact triple on $Q$
\begin{equation}\label{Euler}
0\to\Omega^1_{Q/S}\to p^*(\mathcal{O}_{S}(1)\oplus T_{S}(-1))
\otimes\mathcal{O}_Q(-1)\to\mathcal{O}_Q\to 0,
\end{equation}
we find the $p$-relative dualizing sheaf 
$\omega_{Q/S}:=\det(\Omega^1_{Q/S})$:
\begin{equation}\label{rel dual}
\omega_{Q/S}\simeq\mathcal{O}_Q(-3)\otimes p^*\mathcal{O}_{S}(2).
\end{equation}

Set $\mathcal{E}:=\sigma^*E$. By construction, for each $y\in S$ the 
fiber $Q_y=p^{-1}(y)$ is
a plane such that $l_y=Q_y\cap D_P$ is a line, and, by (\ref{E|Y}),
\begin{equation}\label{triv on l}
\mathcal{E}_{|l_y}\simeq\mathcal{O}_{l_y}^{\oplus2}.
\end{equation}
Furthermore, $\sigma(Q_y)$ is an $\alpha$-plane in $X$, and from 
(\ref{triv on l}) it follows
clearly that 
$h^0(\mathcal{E}_{|Q_y}(-1))=\mathcal{E}^\vee_{|Q_y}(-1))=0$. Hence, in view of condition (a) of the Proposition, the sheaf $\mathcal{E}_{|Q_y}$ is the 
cohomology sheaf of a monad
\begin{equation}\label{eqMonad}
0\to\mathcal{O}_{Q_y}(-1)^{\oplus 
a}\to\mathcal{O}_{Q_y}^{\oplus(2a+2)}\to
\mathcal{O}_{Q_y}(1)^{\oplus a}\to0
\end{equation}
(see \cite[Ch. II, Ex. 3.2.3]{OSS}). This monad immediately implies the 
equalities
\begin{equation}\label{cohomology}
h^1(\mathcal{E}_{|Q_y}(-1))=h^1(\mathcal{E}_{|Q_y}(-2))=a,\ \
h^1(\mathcal{E}_{|Q_y}\otimes\Omega^1_{Q_y})=2a+2,
\end{equation}
$$
h^i(\mathcal{E}_{|Q_y}(-1))=h^i(\mathcal{E}_{|Q_y}(-2))=
h^i(\mathcal{E}_{|Q_y}\otimes\Omega^1_{Q_y})=0,\ \ i\ne1.
$$
Consider the sheaves of $\mathcal{O}_{S}$-modules
\begin{equation}\label{E_i}
E_{-1}:=R^1p_*(\mathcal{E}\otimes\mathcal{O}_Q(-2)\otimes 
p^*\mathcal{O}_{S}(2)),\ \ \
E_0:=R^1p_*(\mathcal{E}\otimes\Omega^1_{Q/S}), \ \ \
E_1:=R^1p_*(\mathcal{E}\otimes\mathcal{O}_Q(-1)).
\end{equation}
The equalities (\ref{cohomology}) together with Cohomology and Base 
Change imply that $E_{-1},\ E_1$ and
$E_0$
are locally free $\mathcal{O}_{S}$-modules, and
$\rk(E_{-1})=\rk(E_1)=a,$ and $\rk(E_0)=2a+2$. Moreover,
\begin{equation}\label{R_i}
R^ip_*(\mathcal{E}\otimes\mathcal{O}_Q(-2))=
R^ip_*(\mathcal{E}\otimes\Omega^1_{Q/S})=R^ip_*(\mathcal{E}\otimes\mathcal{O}_Q(-1))=0
\end{equation}
for $i\ne 1$. Note that $\mathcal{E}^\vee\simeq\mathcal{E}$ as 
$c_1(\mathcal{E})=0$ and
$\rk\mathcal{E}=2$. Furthermore, (\ref{rel dual}) implies that the
nondegenerate pairing ($p$-relative Serre duality)
$R^1p_*(\mathcal{E}\otimes\mathcal{O}_Q(-1))\otimes
R^1p_*(\mathcal{E}^\vee\otimes\mathcal{O}_Q(1)\otimes \omega_{Q/S})\to
R^2p_*\omega_{Q/S}=\mathcal{O}_{S}$ can be rewritten as $E_1\otimes 
E_{-1}\to\mathcal{O}_{S}, $
thus giving an isomorphism
\begin{equation}\label{isom dual}
E_{-1}\simeq E_1^\vee.
\end{equation}
Similarly, since
$\mathcal{E}^\vee\simeq\mathcal{E}$ and $\Omega^1_{Q/S}\simeq 
T_{Q/S}\otimes\omega_{Q/S}$,
$p$-relative Serre duality yields a nondegenerate pairing
$E_0\otimes E_0=R^1p_*(\mathcal{E}\otimes\Omega^1_{Q/S})\otimes
R^1p_*(\mathcal{E}\otimes\Omega^1_{Q/S})=
R^1p_*(\mathcal{E}\otimes\Omega^1_{Q/S})\otimes
R^1p_*(\mathcal{E}^\vee\otimes T_{Q/S}\otimes\omega_{Q/S})
\to R^2p_*\omega_{Q/S}=\mathcal{O}_{S}$.
Therefore $E_0$ is self-dual, i.e. $E_0\simeq E_0^\vee$, and in 
particular $c_1(E_0)=0$.

Now, let $J$ denote the fiber product $Q\times_{S}Q$ with projections
$Q\overset{pr_1}\leftarrow J\overset{pr_2}\to Q$ such that $p\circ 
pr_1=p\circ pr_2$.
Put $F_1\boxtimes F_2:=pr_1^*F_1\otimes pr_2^*F_2$ for sheaves $F_1$ 
and $F_2$ on $Q$,
and consider the standard $\mathcal{O}_J$-resolution of the structure 
sheaf
$\mathcal{O}_{\Delta}$ of the diagonal $\Delta\hookrightarrow J$
\begin{equation}\label{resoln of diag}
0\to\mathcal{O}_Q(-1)\otimes 
p^*\mathcal{O}_{S}(2)\boxtimes\mathcal{O}_Q(-2)\to
{\Omega^1}_{Q/S}(1)\boxtimes\mathcal{O}_Q(-1)\to\mathcal{O}_J\to
\mathcal{O}_{\Delta}\to0.
\end{equation}
Twist this sequence by the sheaf
$(\mathcal{E}\otimes\mathcal{O}_Q(-1))\boxtimes\mathcal{O}_Q(1)$
and apply the functor $R^ipr_{2*}$ to the resulting sequence. In view 
of (\ref{E_i}) and
(\ref{R_i}) we obtain the following monad for $\mathcal{E}$:
\begin{equation}\label{monad1}
0\to p^*E_{-1}\otimes\mathcal{O}_Q(-1)\overset{\lambda}\to 
p^*E_0\overset{\mu}
\to p^*E_1\otimes\mathcal{O}_Q(1)\to0,\ \ \ \ \ \ker(\mu)/{\rm 
im}(\lambda)=\mathcal{E}.
\end{equation}
Put $R:=p^*h$, where $h$ is the class of a line in $S$. Furthermore, 
set
$H:=\sigma^*H_X$, $[\mathbb{P}_\alpha]:=\sigma^*[\mathbb{P}^2_\alpha]$,
$[\mathbb{P}_\beta]:=\sigma^*[\mathbb{P}^2_\beta]$,
where $H_X$ is the class of a hyperplane section of $X$ (via the 
Pl\"ucker embedding),
and respectively, $[\mathbb{P}^2_\alpha]$ and $[\mathbb{P}^2_\beta]$ 
are the classes of an $\alpha$- and
$\beta$-plane. Note that, clearly, 
$\mathcal{O}_Q(H)\simeq\mathcal{O}_Q(1)$.
Thus, taking into account the duality (\ref{isom dual}), we rewrite the 
monad (\ref{monad1}) as
\begin{equation}\label{monad2}
0\to p^*E_1^\vee\otimes\mathcal{O}_Q(-H)\overset{\lambda}\to
p^*E_0\overset{\mu}\to p^*E_1\otimes\mathcal{O}_Q(H)\to0,\ \ \ \ \ \
\ker(\mu)/{\rm im}(\lambda)\simeq\mathcal{E}.
\end{equation}
In particular, it becomes clear that \refeq{monad1} is a relative version of the monad \refeq{eqMonad}.

As a next step, we are going to express all Chern classes of the 
sheaves in
(\ref{monad2}) in terms of $a$. We start by writing down the Chern 
polynomials of the bundles
$p^*E_1\otimes\mathcal{O}_Q(H)$ and 
$p^*E_1^\vee\otimes\mathcal{O}_Q(-H)$ in the form
\begin{equation}\label{Chern1}
c_t(p^*E_1\otimes\mathcal{O}_Q(H))=\prod_{i=1}^a(1+(\delta_i+H)t),\ \ \
c_t(p^*E_1^\vee\otimes\mathcal{O}_Q(-H))=\prod_{i=1}^a(1-(\delta_i+H)t),
\end{equation}
where $\delta_i$ are the Chern roots of the bundle $p^*E_1$. Thus
\begin{equation}\label{c,d}
cR^2=\sum_{i=1}^a\delta_i^2,\ \ dR=\sum_{i=1}^a\delta_i.
\end{equation}
for some $c,d\in\mathbb{Z}$.
Next we invoke the following easily verified relations in $A^\cdot(Q)$:
\begin{equation}\label{rel in A(Q)}
H^4=RH^3=2[pt],\ \ \ R^2H^2=R^2[\mathbb{P}_\alpha]=
RH[\mathbb{P}_\alpha]=H^2[\mathbb{P}_\alpha]=RH[\mathbb{P}_\beta]=H^2[\mathbb{P}_\beta]=[pt],
\end{equation}
$$
[\mathbb{P}_\alpha][\mathbb{P}_\beta]=R^2[\mathbb{P}_\beta]=R^4=R^3H=0,
$$
where $[pt]$ is the class of a point. This, together with (\ref{c,d}), 
gives
\begin{equation}\label{sums}
\sum_{1\le i<j\le a}\delta_i^2\delta_j^2=
\sum_{1\le i<j\le a}(\delta_i^2\delta_j+\delta_i\delta_j^2)H=0,
\sum_{1\le i<j\le a}\delta_i\delta_jH^2=\frac{1}{2}(d^2-c)[pt],
\sum_{1\le i\le a}(\delta_i+\delta_j)H^3=2(a-1)d[pt].
\end{equation}
Note that, since $c_1(E_0)=0$,
\begin{equation}\label{Chern2}
c_t(p^*E_0)=1+bR^2t^2
\end{equation}
for some $b\in\mathbb{Z}$. Furthermore,
\begin{equation}\label{c_t(E)}
c_t(\mathcal{E})=1+a[\mathbb{P}_\alpha]t^2
\end{equation}
by the condition of the Proposition.
Substituting (\ref{Chern2}) and (\ref{c_t(E)}) into the polynomial
$f(t):=c_t(\mathcal{E})c_t(p^*E_1\otimes\mathcal{O}_Q(H))
c_t(p^*E_1^\vee\otimes\mathcal{O}_Q(-H))$, we have
$f(t)=(1+a[\mathbb{P}_\alpha]t^2)\prod_{i=1}^a(1-(\delta_i+H)^2t^2)$.
Expanding $f(t)$ in $t$ and using (\ref{c,d})-(\ref{sums}),
we obtain
\begin{equation}\label{f(t)2}
f(t)=1+(a[\mathbb{P}_\alpha]-cR^2-2dRH-aH^2)t^2+e[pt]t^4,\ \ \
\end{equation}
where
\begin{equation}\label{e}
e=-3c-a(2d+a)+(a-1)(a+4d)+2d^2.
\end{equation}
Next, the monad (\ref{monad2}) implies $f(t)=c_t(p^*E_0)$. A comparison 
of (\ref{f(t)2}) with
(\ref{Chern2}) yields
\begin{equation}\label{c_2}
c_2(\mathcal{E})=a[\mathbb{P}_\alpha]=(b+c)R^2+2dRH+aH^2,
\end{equation}
\begin{equation}\label{c_4}
e=c_4(p^*E_0)=0.
\end{equation}
The relation (\ref{c_4}) is the crucial relation which enables us to 
express the Chern classes
of all sheaves in (\ref{monad2}) just in terms of $a$.
More precisely, (\ref{c_2}) and (\ref{rel in A(Q)}) give
$0=c_2(\mathcal{E})[\mathbb{P}_\beta]=2d+a$, hence
$a=-2d$.
Substituting these latter equalities into (\ref{e}) we get 
$e=-a(a-2)/2-3c$.
Hence $c=-a(a-2)/6$ by (\ref{c_4}). Since $a=-2d$, (\ref{c,d}) and the 
equality $c=-a(a-2)/6$
give
$c_1(E_1)=-a/2,\ \ c_2(E_1)=(d^2-c)/2=a(5a-4)/24$. Substituting this 
into the standard formulas
$e_k:=c_k(p^*E_1\otimes\mathcal{O}_Q(H))=\sum_{i=0}^2\binom{a-i}{k-i}R^iH^{k-i}c_i(E_1),
\ \ 1\le k\le4$, we obtain
\begin{equation}\label{ee_i}
e_1=-aR/2+aH,\ \ e_2=(5a^2/24-a/6)R^2+(a^2-a)(-RH+H^2)/2,
\end{equation}
$$e_3=(5a^3/24-7a^2/12+a/3)R^2H+(-a^3/4+3a^2/4-a/2)RH^2+(a^3/6-a^2/2+a/3)H^3,$$
$$
e_4=(-7a^4/144+43a^3/144-41a^2/72+a/3)[pt].
$$
It remains to write down explicitely $c_2(p^*E_0)$:
(\ref{rel in A(Q)}), (\ref{c_2}) and the relations $a=-2d$, 
$c=-a(a-2)/6$ give
$a=c_2(\mathcal{E})[\mathbb{P}_\alpha]=b+c,$ hence
\begin{equation}\label{c_2(E_0)}
c_2(E_0)=b=(a^2+4a)/6
\end{equation}
by (\ref{Chern2}).

Our next and final step will be to obtain a contradiction by computing 
the Euler characteristic
of the sheaf $\mathcal{E}$ and two different ways. We first compute the 
Todd class
${\rm td}(T_Q)$ of
the bundle $T_Q$. From the exact triple dual to (\ref{Euler}) we find
$c_t(T_{Q/S})=1+(-2R+3H)t+(2R^2-4RH+3H^2)t^2$.
Next, $c_t(T_Q)=c_t(T_{Q/S})c_t(p^*T_S)$.
Hence $c_1(T_Q)=R+3H,\ c_2(T_Q)=-R^2+5RH+3H^2,\ c_3(T_Q)=-3R^2H+9H^2R,\ 
c_4(T_Q)=9[pt].$
Substituting into the formula for
the Todd class of $T_Q$,  ${\rm 
td}(T_Q)=1+\frac{1}{2}c_1+\frac{1}{12}(c_1^2+c_2)+\frac{1}{24}c_1c_2
-\frac{1}{720}(c_1^4-4c_1^2c_2-3c_2^2-c_1c_3+c_4)$, where 
$c_i:=c_i(T_Q)$ (see, e.g.,
\cite[p.432]{H}), we get
\begin{equation}\label{td(T_Q)}
{\rm 
td}(T_Q)=1+\frac{1}{2}R+\frac{3}{2}H+\frac{11}{12}RH+H^2+\frac{1}{12}HR^2+
\frac{3}{4}H^2R+\frac{3}{8}H^3+[pt].
\end{equation}
Next, by the hypotheses of Proposition
$c_1(\mathcal{E})=0,\ c_2(\mathcal{E})=a[\PP_{\alpha}],\ 
c_3(\mathcal{E})=c_4(\mathcal{E})=0$.
Substituting this into the general formula for the Chern character of a 
vector bundle $F$,
$$
{\rm 
ch}(F)=\rk(F)+c_1+(c_1^2-2c_2)/2+(c_1^3-3c_1c_2-3c_3)/6+(c_1^4-4c_1^2c_2+4c_1c_3+2c_2^2-4c_4)/24, \ \
$$
$c_i:=c_i(F)$ (see, e.g., \cite[p.432]{H}), and using
(\ref{td(T_Q)}), we obtain by the Riemann-Roch Theorem for 
$F=\mathcal{E}$
\begin{equation}\label{chi(E)}
\chi(\mathcal{E})=\frac{1}{12}a^2-\frac{23}{12}a+2.
\end{equation}
In a similar way, using (\ref{ee_i}), we obtain
\begin{equation}\label{chi(E1)+chi(E-1)}
\chi(p^*E_1\otimes\mathcal{O}_Q(H))+\chi(p^*E_1^\vee\otimes\mathcal{O}_Q(-H))=
\frac{5}{216}a^4-\frac{29}{216}a^3-\frac{1}{54}a^2+\frac{113}{36}a.
\end{equation}
Next,
in view of (\ref{c_2(E_0)})
and the equality $c_1(E_0)=0$ the Riemann-Roch Theorem for $E_0$
easily gives
\begin{equation}\label{chi(E_0)}
\chi(p^*E_0)=\chi(E_0)=-\frac{1}{6}a^2+\frac{4}{3}a+2.
\end{equation}
Together with (\ref{chi(E)}) and
(\ref{chi(E1)+chi(E-1)}) this yields
$$
\Phi(a):=\chi(p^*E_0)-(\chi(\mathcal{E})+
\chi(p^*E_1\otimes\mathcal{O}_Q(H))+\chi(p^*E_1^\vee\otimes\mathcal{O}_Q(-H)))=
-\frac{5}{216}a(a-2)(a-3)(a-\frac{4}{5}).
$$
The monad (\ref{monad2}) implies now $\Phi(a)=0.$ The only positive 
integer roots of the
polynomial $\Phi(a)$ are $a=2$ and $a=3$.
However, (\ref{chi(E)}) implies $\chi(\mathcal{E})=-\frac{3}{2}$ for 
$a=2$,
and (\ref{chi(E_0)}) implies $\chi(p^*E_0)=\frac{9}{2}$ for $a=3$.
This is a contradiction as the values of $\chi(\mathcal{E})$ and 
$\chi(p^*E_0)$ are integers by
definition.
\end{proof}

We need a last piece of notation.
Consider the flag variety $Fl(k_m-2,k_m+2;V^{n_m})$. Any point
$u=(V^{k_m-2},V^{k_m+2})\in \Fl(k_m-2,k_m+2;V^{n_m})$ determines a 
standard extension
\begin{equation}\label{i_z}
i_{u}:\ X=G(2;4)\hookrightarrow X_m,
\end{equation}
\begin{equation}\label{eq}
W^2\mapsto V^{k_m-2}\oplus W^2\subset V^{k_m+2}\subset
V^{n_m}=V^{k_m-2}\oplus W^4\subset V^{n_m},
\end{equation}
where $W^2\in X=G(2;W^4)$ and an isomorphism $V^{k_m-2}\oplus W^4\simeq 
V^{k_m+2}$ is fixed (clearly $i_{u}$ does not depend on the choice of 
this isomorphism modulo $\Aut(X_m)$). We
clearly have isomorphisms of Chow groups
\begin{equation}\label{isomChow}
i_{u}^*:\ A^2(X_m)\overset{\sim}\to A^2(X),\ \ \
i_{u*}:\ A_2(X)\overset{\sim}\to A_2(X_m),
\end{equation}
and the flag variety $Y_m:=Fl(k_m-1,k_m+1;V^{n_m})$ (respectively, 
$Y:=Fl(1,3;4)$) is the set of lines in $X_m$ (respectively, in $X$).

\vspace{0.3cm}

\begin{theorem}\label{th56}
Let $\displaystyle\GG = \lim_{\to}X_m$ be a twisted ind-Grassmannian. 
Then any vector bundle $\displaystyle\EE=\lim_{\gets}E_m$ on $\GG$ of 
rank 2 is trivial, and hence Conjecture \ref{con1}(iv) holds for vector 
bundles of rank 2.
\end{theorem}

\begin{proof}
Fix $m\ge\max\{m_0,m_1\},$ where $m_0$ and $m_1$ are as in Corollary 
\ref{d=(0,0)} and
Lemma \ref{c_2(E_m)=0}. For $j=1,2$, let $E^{(j)}$ denote the 
restriction of $E_m$ to a
projective plane of type $\mathbb{P}^2_{j,m}$,
 $T^j\simeq\Fl(k_m-j,k_m+3-j,V^{n_m})$ be the variety of planes of the 
form $\mathbb{P}^2_{j,m}$
in $X_m$, and $\Pi^j:=\{\mathbb{P}^2_{j,m}\in T^j|\ 
{E_m}_{|\mathbb{P}^2_{j,m}}$ is properly unstable
(i.e. not semistable)$\}.$ As semistability is an open condition, 
$\Pi^j$ is a closed subset
of $T^j$.

(i) Assume that $c_2(E^{(1)})>0$. Then, since $m\ge m_1$, Lemma 
\ref{c_2(E_m)=0} implies
$c_2(E^{(2)})\le0$.

(i.1) Suppose that $c_2(E^{(2)})=0$. If $\Pi^2\ne T^2$, then for any
$\mathbb{P}^2_{2,m}\in T^2\smallsetminus \Pi^2$ the corresponding 
bundle $E^{(2)}$ is
semistable, hence $E^{(2)}$ is trivial as $c_2(E^{(2)})=0$, see 
\cite[Prop. 2.3,(4)]{DL}. Thus, for a generic point $u\in 
Fl(k_m-2,k_m+2;V^{n_m})$, the bundle $E=i_{u}^*E_m$ on $X=G(2;4)$ satisfies the conditions 
of Proposition \ref{not exist}, which is a contradiction.

We therefore assume $\Pi^2=T^2$. Then for any $\mathbb{P}^2_{2,m}\in 
T^2$ the corresponding bundle $E^{(2)}$ has a maximal destabilizing 
subsheaf $0\to\mathcal{O}_{\mathbb{P}^2_{2,m}}(a)\to E^{(2)}.$ Moreover 
$a>0$. In fact, otherwise the condition $c_2(E^{(2)})=0$ would imply that 
$a=0$ and 
$E^{(2)}/\mathcal{O}_{\mathbb{P}^2_{2,m}}=\mathcal{O}_{\mathbb{P}^2_{2,m}}$, i.e. $E^{(2)}$ would be trivial, in particular 
semistable. Hence
\begin{equation}\label{a,-a}
\mathbf{d}_{E^{(2)}}=(a,-a).
\end{equation}
Since any line in $X_m$ is contained in a plane $\mathbb{P}^2_{2,m}\in 
T^2$,
(\ref{a,-a}) implies $\mathbf{d}_{E_m}=(a,-a)$ with $a>0$ for $m>m_0$, 
contrary to Corollary \ref{d=(0,0)}.

(i.2) Assume $c_2(E^{(2)})<0$. Since $E^{(2)}$ is not stable for any 
$\mathbb{P}^2_{2,m}\in T^2$, its maximal destabilizing subsheaf 
$0\to\mathcal{O}_{\mathbb{P}^2_{2,m}}(a)\to E^{(2)}$ clearly satisfies the 
condition $a>0$, i.e. $E^{(2)}$ is properly unstable, hence $\Pi^2=T^2$. 
Then we again obtain a contradiction as above.

(ii) Now we assume that $c_2(E^{(2)})>0$. Then, replacing $E^{(2)}$ by 
$E^{(1)}$ and vice versa, we arrive to a contradiction by the same 
argument as in case (i).

(iii) We must therefore assume $c_2(E^{(1)})=c_2(E^{(2)})=0$. Set 
$D(E_m):=\{l\in Y_m|~\mathbf{d}_{E_m}(l)\ne(0,0)\}$ and $D(E):=\{l\in 
Y|~\mathbf{d}_E(l)\ne(0,0)\}$. By Corollary \ref{d=(0,0)}, 
$\mathbf{d}_{E_m}=(0,0),$ hence $\mathbf{d}_E=(0,0)$ for a generic embedding $i_u:X\hookrightarrow X_m$. Then by deformation theory \cite{B}, $D(E_m)$ (respectively, 
$D(E)$) is an effective divisor on $Y_m$ (respectively, on $Y$). Hence, 
$\mathcal{O}_Y(D(E))=p_1^*\mathcal{O}_{Y^1}(a) \otimes 
p_2^*\mathcal{O}_{Y^2}(b)$ for some $a,b\ge0$, where $p_1$, $p_2$ are as in diagram 
\refeq{eqDiag}. Note that each fiber of $p_1$ (respectively, of $p_2$) is a plane 
$\tilde{\mathbb{P}}^2_{\alpha}$ dual to some $\alpha$-plane 
$\mathbb{P}^2_{\alpha}$ (respectively, a plane $\tilde{\mathbb{P}}^2_{\beta}$ dual to 
some $\beta$-plane $\mathbb{P}^2_{\beta}$). Thus, setting
$D(E_{|\mathbb{P}^2_{\alpha}}):=\{l\in\tilde{\mathbb{P}}^2_{\alpha}|~\mathbf{d}_E(l)\ne(0,0)\}$, 
$D(E_{|\mathbb{P}^2_{\beta}}):=\{l\in\tilde{\mathbb{P}}^2_{\beta}|~\mathbf{d}_E(l)\ne(0,0)\}$,
we obtain
$\mathcal{O}_{\tilde{\mathbb{P}}^2_{\alpha}}(D(E_{|\mathbb{P}^2_{\alpha}}))=
\mathcal{O}_Y(D(E))_{|\tilde{\mathbb{P}}^2_{\alpha}}=
\mathcal{O}_{\tilde{\mathbb{P}}^2_{\alpha}}(b),\ \ \
\mathcal{O}_{\tilde{\mathbb{P}}^2_{\beta}}(D(E_{|\mathbb{P}^2_{\beta}}))=
\mathcal{O}_Y(D(E))_{|\tilde{\mathbb{P}}^2_{\beta}}=
\mathcal{O}_{\tilde{\mathbb{P}}^2_{\beta}}(a).$ Now if
$E_{|\mathbb{P}^2_{\alpha}}$ is semistable, a theorem of
Barth \cite[Ch. II, Theorem 2.2.3]{OSS} implies that
$D(E_{|\mathbb{P}^2_{\alpha}})$ is a divisor of degree
$c_2(E_{|\mathbb{P}^2_{\alpha}})=a$ on
$\mathbb{P}^2_{\alpha}$.
Hence $a=c_2(E^{(1)})=0$ for a semistable $E_{|\mathbb{P}^2_{\alpha}}$.
If $E_{|\mathbb{P}^2_{\alpha}}$ is not semistable, it is unstable and 
the equality
$\mathbf{d}_E(l)=(0,0)$ yields 
$\mathbf{d}_{E_{|\mathbb{P}^2_{\alpha}}}=(0,0)$. Then the
maximal destabilizing subsheaf of $E_{|\mathbb{P}^2_{\alpha}}$ is 
isomorphic to
$\mathcal{O}_{\mathbb{P}^2_{\alpha}}$ and, since 
$c_2(E_{|\mathbb{P}^2_{\alpha}})=0,$
we obtain an exact triple
$0\to\mathcal{O}_{\mathbb{P}^2_{\alpha}}\to
E_{|\mathbb{P}^2_{\alpha}}\to \mathcal{O}_{\mathbb{P}^2_{\alpha}}\to 
0$,
so that
$E_{|\mathbb{P}^2_{\alpha}}\simeq\mathcal{O}_{\mathbb{P}^2_{\alpha}}^{\oplus2}$ 
is semistable,
a contradiction. This shows that $a=0$ whenever 
$c_2(E^{(1)})=c_2(E^{(2)})=0$. Similarly, $b=0$.
Therefore $D(E_m)=\emptyset$, and
Proposition \ref{prop31} implies that $E_m$ is trivial. Therefore 
$\mathbf{E}$ is trivial as
well.
\end{proof}

In \cite{DP} Conjecture \ref{con1} (iv) was proved
not only when $\mathbf{X}$ is a twisted projective ind-space,
but also for finite rank bundles on special twisted ind-Grassmannians 
defined through certain
homogeneous embeddings $\phi_m$. These include embeddings of the form
\[
G(k;n)\to G(ka;nb)
\]
\[
V^k\subset V\mapsto V^k\otimes W^a\subset V\otimes W^b,
\]
where $W^a\subset W^b$ is a fixed pair of finite-dimensional spaces 
with $a>b$, or of the form
\[
G(k;n)\to G\left(\frac{k(k+1)}{2};n^2\right)
\]
\[
V^k\subset V\mapsto S^2(V^k)\subset V\otimes V.
\]
More precisely, Conjecture \ref{con1} (iv) was proved in \cite{DP} for 
twisted
ind-Grassmannians whose defining embeddings are homogeneous embeddings 
satisfying some specific
numerical conditions relating the degrees $\deg\phi_m$ with the pairs 
of integers $(k_m,n_m)$.
There are many twisted ind-Grassmannians for which those conditions are 
not satisfied. For
instance, this applies to the ind-Grassmannians defined by iterating 
each of the following
embeddings:
\begin{eqnarray*}
G(k;n)\to G\left(\frac{k(k+1)}{2};\frac{n(n+1)}{2}\right)\\
V^k\subset V\mapsto S^2(V^k)\subset S^2(V), \\
G(k;n)\to G\left(\frac{k(k-1)}{2};\frac{n(n-1)}{2}\right)\\
V^k\subset V\mapsto  \Lambda^2(V^k)\subset \Lambda^2(V).
\end{eqnarray*}
Therefore the resulting ind-Grassmannians $\mathbf G(k,n,S^2)$ and
$\mathbf G (k,n,\Lambda^2)$ are examples of twisted ind-Grassmannians 
for which \refth{th56}
is new.

\end{document}